\documentclass[12pt]{article}
\usepackage{amssymb,amscd,amsmath,amsthm}
\usepackage{latexsym,amstext}
\usepackage{latexsym,amstext}
\usepackage{color}
\usepackage{graphicx}
\usepackage{epstopdf}
\usepackage{cite}

\begin{document}

\title{\bf   A note on linear differential equations with variable coefficients }

\author{M. Gadella$^1$, L.P. Lara$^{2,3}$\\ \\
$^1$ Departamento de F\'{\i}sica Te\'orica, At\'omica y Optica  and IMUVA, \\
Universidad de Va\-lladolid, 47011 Valladolid, Spain\\ 
$^2$ Instituto de F\'isica Rosario, CONICET-UNR, \\ 
Bv. 27 de Febrero, S2000EKF Rosario, Santa Fe,  Argentina.\\
$^3$ Departamento de Sistemas, Universidad del Centro Educativo\\
Latinoamericano, Av. Pellegrini 1332, S2000 Rosario, Argentina
}

\maketitle

\begin{abstract}

In this manuscript, we deal with some particular type of homogeneous first order linear systems with variable coefficients, in which we provide qualitative properties of the solution. When the coefficients of the indeterminate functions are periodic with the same period, $T$, we obtain a simple method so as to obtain the Floquet coefficients. We give a new interpretation of the averages on the interval $(0,T)$ of the matrix of the coefficients. We discuss some examples and made some comparison with previous results.

\end{abstract}

\section{Introduction and Motivation}

Ordinary differential equations or systems with variable coefficients appear quite often in Physics or Engineering. Except for a bunch of well studied cases, to find exact solutions for these systems is not easy in general, even if the equations are linear. See \cite{HAL,CHI,LEF,PER,YS,VER}. Of particular interest are the linear two dimensional systems, which in particular include linear equations of the form 

\begin{equation}\label{1}
a(t)\, y''(t) + b(t)\, y'(t) + c(t)\, y(t) =0\,,
\end{equation}
which is equivalent to the system ($z=y'$)

\begin{equation}\label{2}
\left( \begin{array}{c} y' \\ [2ex] z'  \end{array} \right)  \left( \begin{array}{cc} 1 & 0 \\ [2ex] -\frac{b}{a} & - \frac{c}{a} \end{array} \right)  \left( \begin{array}{c} z \\ [2ex] y \end{array} \right)\,.
\end{equation}

Equations of type \eqref{1} are ubiquitous in classical and, particularly, quantum physics. Equations and systems of equations with periodic coefficients are of special interest as they describe the behaviour of periodic physical systems.

This note intends to be a contribution to the analysis of two dimensional linear equations of the form

\begin{equation}\label{3}
\dot{\mathbf x}(t) = A(t) \mathbf x(t)\,,
\end{equation}
where the dot means the derivative with respect to the variable $t$.  We look for solutions with initial condition $\mathbf x(t_0)= \mathbf x_0$, $\mathbf x_0 \in \mathbb R^2$. The entries of $A(t)$, $\{a_{ij}(t)\}$, are real analytic functions of $t$ on a neighbourhood of $t_0$. Although the restriction to $2 \times 2$ $A(t)$ matrices may look important, this includes two  cases of particular interest. The former is the one dimensional, in the coordinate variable, two dimensional in phase space, classical dynamical systems. The other is the one dimensional Schr\"odinger equation for quantum systems. Both are interesting as laboratories to study the behaviour of dynamical systems and in general their dynamical equations involve less difficulties on their resolution than those for higher dimensions. 

In spite of their apparent simplicity, equations of the form \eqref{3} cannot be solved in general. Therefore, some restrictions are in order even to find qualitative properties of their solutions. In this note, we propose some restrictions based in commutativity conditions, which permit to obtain explicitly the form of the fundamental matrix, $\Phi(t)$, of \eqref{3}.  Limit values of the norm of $\Phi(t)$ in terms of the variable $t$ can be obtained, which eventually would permit a discussion on the stability of solutions.  

Precisely, Lyapunov stability of solutions is the object of study of an article \cite{GRAU}, which in part has motivated the present note. This article stresses the fact that, in spite of their interest, not much investigation has been published on systems like that in \eqref{3}, even for the case in which the coefficients of the matrix $A(t)$ were periodic with the same period $T$. We have added a comment of such systems in our Section 3. Under our hypothesis, we have obtained explicit expressions of the Floquet coefficients in terms of averages of the entries of $A(t)$. 

Finally, we give some examples in which we obtain explicit solutions to system \eqref{3}. These examples show that our commutation hypothesis are far from trivial and also that can be find exact solutions for models which cannot be obtained using symbolic calculus methods.

\section{Two dimensional systems}

Let us go back to system \eqref{3}, where we assume that $A(t)$ is a $2 \times 2$ matrix. Then, let us define the following matrix:

\begin{equation}\label{4}
D(t):= \int_{t_0}^t A(s)\, ds\,,
\end{equation}
where the integral means the matrix obtained by integration of the entries of $A(s)$. Then, obviously, $\dot D(t)=A(t)$, i.e., the derivative of $D(t)$ with respect to $t$ coincides with $A(t)$. 

We depart from the following Ansatz: The matrices $A(t)$ and $D(t)$ commute, $A(t)D(t)= D(t)A(t)$. A sufficient condition to verify this Ansatz is that $D(t)$ has the following form:

\begin{equation}\label{5}
D(t) = \left(\begin{array}{cc} f(t) & g(t) \\[2ex] \alpha g(t) & f(t) + \beta g(t) \end{array}\right)\,,
\end{equation}
where $f(t)$ and $g(t)$ are real differentiable functions and $\alpha$ and $\beta$ are real numbers. In terms of the entries of the matrix $A(t)$, functions $f(t)$ and $g(t)$ have the following form:

\begin{equation}\label{6}
f(t) = \int_{t_0}^t a_{11}(s)\, ds \,, \qquad g(t) = \int_{t_0}^t a_{12}(s)\, ds\,.
\end{equation}

Next, we assume a second Ansatz: The matrices $A(t)$ and its first derivative $\dot A(t)$, which is the matrix for which the entries are the derivatives of the entries of $A(t)$, commute: $A(t)\dot A(t) = \dot A(t) A(t)$. Exactly as in \eqref{5}, a sufficient condition for this commutativity is that

\begin{equation}\label{7}
A(t) = \left(\begin{array}{cc} a_{11}(t) & a_{12}(t) \\ [2ex]  \alpha\, a_{12}(t)  & a_{11}(t) + \beta\, a_{12}(t)   \end{array}\right)
\end{equation}

Under the above hypothesis, a fundamental matrix, $\Phi(t)$, of \eqref{3} is given by

\begin{equation}\label{8}
\Phi(t) = e^{D(t)}\,,
\end{equation}
so that the solution with the given initial value is

\begin{equation}\label{9}
\mathbf x(t) = \Phi(t) \, \mathbf x_0\,.
\end{equation}

Let us determine the explicit form of the fundamental matrix \eqref{8}. First of all, we write $D(t)$ as

\begin{equation}\label{10}
D(t) = \left(f(t) + \frac \beta 2\, g(t)   \right)I + g(t) \left(\begin{array}{cc} -\frac\beta 2  & 1 \\[2ex] \alpha & \frac \beta 2 \end{array} \right)  \,.
\end{equation}

Since one of the above summands is a multiple of the identity matrix, the exponential of the sum is the product of the exponentials:

\begin{equation}\label{11}
\Phi(t) = \exp \left[  f(t) + \frac \beta 2\, g(t)  \right] \exp\left[  g(t) \left(\begin{array}{cc} -\frac\beta 2  & 1 \\[2ex] \alpha & \frac \beta 2 \end{array} \right)  \right]\,.
\end{equation}

Let us call $S(t)$ to the second term in \eqref{10}. In order to explicitly obtain its exponential, second factor in \eqref{11}, we should diagonalize it. Its eigenvalues are

\begin{equation}\label{12}
\nu_\pm = \pm \gamma g(t)\,, \quad \gamma = + \frac 12 \sqrt{4\alpha^2 + \beta^2}\,.
\end{equation}

For $\gamma \ne 0$, $\alpha \ne 0$, their corresponding eigenvectors are

\begin{equation}\label{13}
v_\pm = \frac 1 \alpha \left( -\frac \beta 2 \pm \gamma, 1 \right)\,.
\end{equation}

If we call $S_d(t)$ to the diagonal form of $S$, then, $S_d= Q^{-1} S Q$, with 

\begin{equation}\label{14}
Q= \left(\begin{array}{cc} \frac 1 \alpha \left(- \frac\beta 2 + \gamma  \right) & \frac 1 \alpha \left(- \frac\beta 2 - \gamma  \right) \\ [2ex] 1  & 1 \end{array}\right)\,.
\end{equation}

Thus, 

\begin{equation}\label{15}
e^{S(t)} = Q\, e^{S_d(t)} Q^{-1} = [2\gamma \cosh (\gamma g(t)) - \beta \sinh(\gamma g(t))] I  + 2 \sinh(\gamma g(t)) \left(\begin{array}{cc} 0 & 1 \\ [2ex] \alpha & 0   \end{array}\right)\,,
\end{equation}
where $I$ is the $2\times 2$ identity matrix. With \eqref{11}, this gives the form of $\Phi(t)$ for $\gamma \ne 0$. In order to find the fundamental matrix for $\gamma=0$, we just take the limit $\gamma\longmapsto 0$, so as to obtain

\begin{equation}\label{16}
\Phi(t) = \exp \left[  f(t) + \frac \beta 2\, g(t)  \right]  \left(\begin{array}{cc}  1- \frac \beta 2\, g(t)  & g(t)  \\ [2ex] - \frac{\beta^2}{4}\, g(t)  & 1+ \frac \beta 2\, g(t)  \end{array}\right)\,,
\end{equation}
where we have still kept $\alpha\ne 0$. If $\alpha =0$ and $\gamma \ne 0$, the eigenvalues of $S(t)$ are just those given in \eqref{12}, with respective eigenvectors 

\begin{equation}\label{17}
v_+ = (1,0)\,, \qquad v_-=\left( \frac1\beta,1 \right)\,.
\end{equation}

Here, the fundamental matrix is obtained by choosing $\alpha= 0$ in \eqref{15} and multiplying the result for the first factor in \eqref{11}. 

\subsection{Some properties}

Here, we make some comments. The three first refer to properties than come after the assignments of some values for $\gamma$ and are straightforward consequences of the form of the fundamental matrix $\Phi(t)$. The forth one is just a particular case.

\begin{itemize}

\item{If $\gamma>0$, we have that

\medskip

1.- If $\beta>0$, $f(t)$ is bounded and $\lim_{t\to\infty}g(t)=\infty$, then, $\lim_{t\to\infty}  ||\Phi(t)||\longmapsto \infty$.

\medskip

2.- If $\lim_{t\to\infty}f(t)=\infty$ and $g(t)$ is bounded, then also, $\lim_{t\to\infty}  ||\Phi(t)||\longmapsto \infty$.

\medskip

3.- If $\lim_{t\to\infty}f(t)=-\infty$ and $g(t)$ is bounded, then, $\lim_{t\to\infty}  ||\Phi(t)||\longmapsto 0$.

}

\item{If $\gamma$ is imaginary,

\medskip

1.- If $ f(t) + \frac \beta 2\, g(t)=0$, then $||\Phi(t)||$ is periodic. 

\medskip

2.- If $f(t) + \frac \beta 2\, g(t)$ is periodic, then $||\Phi(t)||$ is either periodic or quasiperiodic. 
}

\item{If $\gamma=0$,

\medskip

1.- If both, $f(t)$ and $g(t)$ were bounded, then $||\Phi(t)||$ is bounded. 

\medskip

2.- If $\lim_{t\to\infty}f(t)= 0$ and $\lim_{t\to\infty}g(t)=0$, then $\lim_{t\to\infty} ||\Phi(t)||=0$. 

\medskip

3.- If $\lim_{t\to\infty}f(t)= \infty$, then $\lim_{t\to\infty} ||\Phi(t)||=\infty$. 

\medskip

4.- If $\lim_{t\to\infty}f(t)= -\infty$, then $\lim_{t\to\infty} ||\Phi(t)||=0$. 

\medskip
}

\item{{\it A particular case}

\medskip

Let us consider the equation:

\begin{equation}\label{18}
\dot {\mathbf  x}(t) = A(t) \mathbf x(t) + B(\mathbf x(t),t)\,,
\end{equation}
where $A(t)$ has the form \eqref{5} and $B(\mathbf x(t),t)$ is a $2 \times 2$ matrix continuous on its arguments. We are considering the initial value $\mathbf x(0)= \mathbf x_0$. If $\Phi(t)$ is the fundamental matrix corresponding to its homogeneous part, we may proceed with the change of indeterminate:

\begin{equation}\label{19}
\mathbf y(t):= \Phi^{-1}(t) \mathbf x(t)  \Longleftrightarrow \mathbf x(t) = \Phi(t) \mathbf y(t) \,.
\end{equation}

Deriving the second form of \eqref{19} with respect to $t$ and taking into account that $\dot\Phi(t)= A(t)\Phi(t)$, we readily obtain that

\begin{equation}\label{20}
\dot{\mathbf y}(t) = \Phi^{-1}(t) B(\Phi(t) \mathbf y(t),t)\,,
\end{equation}
where the given initial condition is equivalent to $\mathbf y(0) = \Phi^{-1}(0) \mathbf x_0$. Thus, in order to find the solution $\mathbf x(t)$ of \eqref{18}, it is sufficient the integration of system \eqref{20}. 
}

\end{itemize}

\section{Systems with periodic coefficients}

In the present Section, we make some comments on systems with periodic coefficients. Typically, these are systems of the form  \eqref{1}, where the matrix $A(t)$ is $n \times n$, real and periodic, all its entries are periodic, with the same period $T$, $A(t+T)=A(t)$. The fundamental matrix may be written as

\begin{equation}\label{21}
\Phi(t) = P(t) \, e^{Bt}\,.
\end{equation}

Here, $B$ is a constant matrix.  We know that there is a non-singular matrix $C$ with

\begin{equation}\label{22}
\Phi(t+T) = \Phi(t)\, C\,.
\end{equation}

A necessary and sufficient condition for $P(t)$ be periodic is that $C$ be of the form $C=e^{BT}$. The proof is trivial.

As is well known, the {\it caracteristic multipliers} are the eigenvalues of $C$, while the {\it characteristic exponents}, also called the {\it Floquet exponents} are the eigenvalues of $B$. We denote characteristic multipliers and exponents by $\rho$ and $\lambda$, respectively. Note the relation between them:

\begin{equation}\label{23}
\rho=e^{\lambda T}\,.
\end{equation}

The imaginary part of $\lambda$ is defined modulo $2\pi$. Then, we define $\mathbf y(t)$ as (compare with \eqref{19}):

\begin{equation}\label{24}
\mathbf y(t):= P(t) \mathbf x(t)  \Longleftrightarrow \mathbf x(t) = P(t) \mathbf y(t)\,.
\end{equation}

Then, \eqref{1} and \eqref{24} give

\begin{equation}\label{25}
\dot P(t) \mathbf y(t) + P(t) \dot {\mathbf y}(t) = A(t) P(t) \mathbf y(t)\,,
\end{equation}
so that

\begin{equation}\label{26}
\dot {\mathbf y}(t) = P^{-1}(t) [A(t) P(t) - \dot P(t)] \mathbf y(y)\,.
\end{equation}
 
Then, take the derivative with respect to $t$ of the expression $P(t) = \Phi(t) e^{-Bt}$ and do some algebra so as to obtain:

\begin{equation}\label{27}
\dot P(t) = A(t) P(t) -P(t) B\,.
\end{equation}

Using \eqref{25} in \eqref{26}, we finally obtain,

\begin{equation}\label{28}
\dot {\mathbf y}(t) = B \mathbf y(t)\,,
\end{equation}
which is a system with constant coefficients, for which the solution procedure is well known. From the form of the Floquet coefficients, we determine the stability of the solution $\mathbf y(t)= \mathbf 0$ and, hence, of all solutions.  

Nevertheless, the matrix $B$ may rarely be determined analytically, so that numerical methods are in order so as to obtain the Floquet exponents. Recall that for the Floquet exponents and for the characteristic multipliers, the following relations hold, respectively ($n$ is the total number of parameters coinciding with the dimension of the matrices $B$, etc):

\begin{equation}\label{29}
\sum_{k=1}^n = \frac 1T \int_0^T {\rm tr}\, A(t)\, dt\,,
\end{equation}
modulo $(2\pi i)/T$ and

\begin{equation}\label{30}
\prod_{k=1}^n \rho_k = \exp \left[ \int_0^T {\rm tr}\, A(t)\, dt \right]\,.
\end{equation}

\section{Operations with commuting matrices.}

Let us go back to a periodic $2 \times 2$ $A(t)$ with period $T$ and initial condition $\mathbf x(t_0)= \mathbf x_0$, $\mathbf x_0 \in \mathbb R^2$. From \eqref{22}, we have that

\begin{equation}\label{31}
C= \Phi^{-1}(t_0) \Phi(t_0+T)\,.
\end{equation}

Go back to \eqref{8}. This gives $\Phi^{-1}(t)=e^{-D(t)}$. According to \eqref{5} and \eqref{6}, $D(t_0)=I$, the identity matrix. We always may choose $t_0=0$ by simplicity, which gives

\begin{equation}\label{32}
C = \Phi(T)\,,
\end{equation}
so that \eqref{32}, \eqref{6} and $C=e^{BT}$ give

\begin{equation}\label{33}
e^{D(T)} = e^{BT}  \Longleftrightarrow  B= \frac 1T \, D(T)\,.
\end{equation}

This means that the matrix $B$ represents the averages on time of $A(t)$. Since the Floquet exponents are the eigenvalues of the matrix $B$, equation \eqref{33} shows two interesting advantages. First, we may determine rather easily the Floquet exponents as eigenvalues of $B$. In addition, we find an interesting interpretation to the entries of $A(t)$. The Floquet exponents may be written in terms of averages of the entries of $A(t)$ as

\begin{equation}\label{34}
\lambda_\pm = \overline{a_{11}}  +  \frac \beta 2 \, \overline{a_{12}} \pm \gamma \overline{a_{12}}\,,
\end{equation}
where,

\begin{equation}\label{35}
\overline{a_{11}}  = \frac 1T \, f(T)\,, \qquad \overline{a_{12}} = \frac 1T \, g(T)\,,
\end{equation}
are the averages of $a_{11}(t)$ and $a_{12}(t)$ on the interval $(0,T)$, since the functions $f(t)$ and $g(t)$ are defined as in \eqref{6} and we have chosen $t_0=0$. 

Some particular cases of interest:

\begin{itemize}

\item{Let $\gamma$ be purely imaginary and $\overline{a_{12}}\ne 0$. Then, if $\overline{a_{11}}  +  \frac \beta 2 \, \overline{a_{12}}$ is positive, there are divergent solutions. However, if it were negative there are solutions that oscillate to zero as $t \longmapsto \infty$.   }

\item{ If $\overline{a_{11}}  +  \frac \beta 2 \, \overline{a_{12}} =0$ all solutions are bounded.}

\item{ If $\gamma$ were real, there is no oscillatory solutions. }

\item{ If $\alpha =0$, the Floquet exponents are

\begin{equation}\label{36}
\lambda_+ = \overline{a_{11}}  +  \frac \beta 2 \, \overline{a_{12}}\,, \qquad \lambda_- =  \overline{a_{11}} \,.
\end{equation}
}

\end{itemize}

In summary, we have proven that {\it the Floquet exponents are the eigenvalues of the matrix $\overline{A(t)}$}, which is the matrix with entries equal to the averages of the entries of $A(t)$. 

\section{Some examples}

The first example comes from a system from \cite{GRAU}, where the authors consider a two dimensional system \eqref{1} with matrix $A(t)$ given by

\begin{equation}\label{37}
A(t) = \left( \begin{array}{cc} \sigma_0 -a(t) & - \sigma_2 \, a(t) \\ [2ex] \sigma_1 \, a(t) & \sigma_0 + a(t) \end{array} \right)\,,
\end{equation}
where $\sigma_i$, $i=0,1,2,$, are real numbers. The function $a(t)$ is real, continuous and $T$-periodic. This matrix has the form \eqref{7} with

\begin{equation}\label{38}
\alpha = - \frac{\sigma_1}{\sigma_2}\,, \qquad \beta = - \frac{2}{\sigma_2} \,, \qquad a_{11}= \sigma_0 -a\,, \qquad a_{12} = - \sigma_2 \, a\,,
\end{equation}
so that the system \eqref{1} with matrix $A(t)$ given by \eqref{35} has a solution of the form \eqref{14}. Floquet exponents are given by \eqref{34} and have here the following explicit form 

\begin{equation}\label{39}
\lambda_\pm = \sigma_0 \pm \frac 1T \, \sqrt{1- \sigma_1 \, \sigma_2} \int_0^t a(s)\, ds \,.
\end{equation}

This shows that the zero solution of \eqref{1} with \eqref{37} is asymptotically estable if $\sigma_0<0$ and $\sigma_1 \, \sigma_2 >1$. This is a condition less restrictive that the one given in \cite{GRAU}, which is $\sigma_0<0$, $\sigma_1<0$, $\sigma_2<0$ and $\sigma_1 \, \sigma_2 >1$.

Particular solutions are determined by initial conditions. For instance, if the initial condition were $\mathbf x(0) = (0,1)$, the Floquet exponents would be  $\lambda_\pm = -1 \pm \frac 12\, i$ and the particular solution has the following components:

\begin{eqnarray}\label{40}
x_1(t) &=& - e^{-t} \, \sin \left( \frac 12 \, t + \frac 12 \, \sin t \, \cos t \right) \,, \nonumber \\ [2ex] 
x_2 (t)  &=& e^{-t} \, \left[ \cos \left( \frac 12 \, t + \frac 12 \, \sin t\, \cos t  \right) - \sin \left(  \frac 12 \, t + \frac 12 \, \sin t\, \cos t \right) \right] \,.
\end{eqnarray}

The solutions are defined on the range of values of $t$, at which $a(t)$ be continuous. 

\medskip

As a second example, let us take \eqref{1} with $A(t)$ given as in \eqref{5}. If $\mathbf x(t)$ is a given solution, let us use the notation $\mathbf x(t) = (x(t), y(t))$. Then, \eqref{3} with \eqref{7} is equivalent to a second order equation

\begin{equation}\label{41}
\ddot x(t) - (2a_{11} + s(t)) \dot x(t) +(a_{11}^2 - \alpha a_{12} - \dot a_{11} + a_{11}\, s(t) ) x(t)=0\,,
\end{equation}
with 

\begin{equation}\label{42}
s(t):= \beta a_{12} + \frac{\dot a_{12}}{a_{12}}\,.
\end{equation}

Many examples have shown that \eqref{41} cannot be, in general, solved through a symbolic calculus program. However, for given $a_{11}(t)$ and $a_{12}(t)$ integrable, our method allows to find an explicit solution. Take as an example

\begin{equation}\label{43}
a_{11}(t)= - \frac{1}{2}\, s(t)\,, \qquad a_{12}(t) = e^{t^2}\,.
\end{equation}

Using these parameters, equation \eqref{43} cannot be solved by symbolic calculus. However, if we use the method described in the present manuscript we obtain an explicit solution as:

\begin{equation}\label{44}
x(t) = \frac{1}{2\gamma} \, \exp \left\{2 \gamma \, \cosh \left( \frac{\gamma}{2} \, \sqrt \pi \, {\rm erf}_i (t) \right) - \beta\, \sinh \left( \frac{\gamma}{2} \, \sqrt \pi \, {\rm erf}_i (t) \right)   \right\}\, 
\end{equation}
where ${\rm erf}_i (t) =-i \,{\rm erf}(it)$, ${\rm erf}(t)$ being the error function. 

\medskip

Let us go back to equation \eqref{41} and let us fix $a_{11}(t) = -\beta\, a_{12}(t)$. We obtain the so call parametric oscillator, with equation

\begin{equation}\label{45}
\ddot x(t) = \left( \beta\, a_{12}(t) - \frac{\dot a_{12}(t)}{a_{12}(t)} \right)\, \dot x(t) - \alpha\, a_{12}^2(t) \, x(t) =0\,.
\end{equation}

Let us take some particular values such as

\begin{equation}\label{46}
\alpha = -1\,, \qquad \beta= \frac \nu \omega \,, \qquad a_{12}(t) = \omega\,,
\end{equation}
where $\nu$ and $\omega$ are fixed positive constants, we obtain the damped linear oscillator as

\begin{equation}\label{47}
\ddot x(t) + \nu\, \dot x(t) + \omega^2\, x(t) =0\,,
\end{equation}
for which the method of solution is well known. 

\section{Concluding remarks}

This is a contribution to the methods of resolution of one dimensional dynamical systems both classical and quantum as they dynamical equations depends on variable coefficients. We have obtained general solutions under some nontrivial conditions based on commutation relations given by the matrix of the coefficients, $A(t)$, with other matrices derived from them such as the derivative or the primitive of $A(t)$. These matrices have entries which are derivatives or primitive of the coefficients of $A(t)$. These type of restrictions, like some others, are somehow necessary as the problem of the solvability of equation \eqref{3} is not easy in general, even for periodic systems. With our hypothesis, we may solve some one dimensional systems which are not solvable using symbolic calculus. 

This note is a part of a research carried out by the authors  on the subject. In \cite{GL}, we have proposed a method to obtain local analytic approximate solutions of ordinary differential equations with variable coefficients, or even some nonlinear equations, inspired in the Lyapunov method, where instead of polynomial approximations, we use truncated Fourier series with variable coefficients as approximate solutions. In \cite{GL1},  we have proposed a modification of a method based on Fourier analysis to obtain the Floquet characteristic exponents for periodic homogeneous linear systems,
which shows a high precision.

\section* {Acknowledgements}

The paper has been parcially supported by MCIN of Spain with funding from the European Union NextGenerationEU
(PRTRC17.I1) and the PID2020-113406GB-IO project by MCIN.

\end{document}